\documentclass[11pt,a4paper]{amsart}
\usepackage{url,amssymb,latexsym}
\usepackage{mymath,mythm}  
\usepackage{vmargin}
\newcommand{\zerovec}{\ensuremath{\mathbf{0}}}
\newcommand{\onevec}{\ensuremath{\mathbf{1}}}
\newcommand{\cvec}{\ensuremath{\mathbf{c}}}

\newcommand{\nvec}{\ensuremath{\mathbf{n}}}

\newcommand{\xvec}{\ensuremath{\mathbf{x}}}
\newcommand{\xhat}{\ensuremath{\mathbf{\widehat{x}}}}
\newcommand{\hmatrix}{\ensuremath{\mathbf{H}}}
\newcommand{\thetahat}{\ensuremath{\mathbf{\widehat{\theta}}}}
\newcommand{\ftilde}{\ensuremath{\widetilde{f}}}
\DeclareMathOperator{\Dif}{D}
\DeclareMathOperator{\Res}{Res}
\newcommand{\s}[1]{\ensuremath{\mathcal{#1}}}
\newcommand{\crit}[1]{\ensuremath{\textsc{Crit(#1)}}}
\renewcommand{\contrib}[1]{\ensuremath{\textsc{Contrib(#1)}}}
\newcommand{\exsym}{\hfill\ensuremath{\lhd}}
\begin{document}
\frenchspacing
\title[A new diagonal method]
  {A new method for computing asymptotics of diagonal coefficients 
of multivariate generating functions}
\author{Alexander Raichev}
\address{Department of Computer Science \\ University of Auckland \\
  Private Bag 92019 \\ Auckland \\ New Zealand}
\email{raichev@cs.auckland.ac.nz}
\author{Mark C. Wilson}
\address{Department of Computer Science \\ University of Auckland \\
  Private Bag 92019 \\ Auckland \\ New Zealand}
\email{mcw@cs.auckland.ac.nz}
\begin{abstract}
Let $\sum_{\nvec\in\nats^d} f_\nvec \xvec^\nvec$ be a multivariate generating
function that converges in a neighborhood of the origin of $\comps^d$.
We present a new, multivariate method for computing the asymptotics of the
diagonal coefficients $f_{a_1n,\ldots,a_dn}$ and show its superiority over
the standard, univariate diagonal method.
\end{abstract}
\maketitle

\section{Introduction}

Let $F(\xvec) = \sum_{\nvec\in\nats^d} f_\nvec \xvec^\nvec$ be a complex 
power series that converges in a neighborhood of the origin 
but not on all of $\comps^d$.
Here $\xvec = (x_1,\ldots,x_d)$, $\nvec = (n_1,\ldots,n_d)$, and
$f_\nvec \xvec^\nvec = f_{n_1,\ldots,n_d} x_1^{n_1}\cdots x_d^{n_d}$.

We wish to compute asymptotics for the diagonal coefficients
$f_{a_1n,\ldots,a_dn}$ for fixed positive integers $a_1,\ldots,a_d$,
a task often useful in enumerative combinatorics. For simplicity of presentation, we 
suppose that $f_\nvec\ge 0$ and that $F$ is rational, although much greater generality is possible.
Thus far the only general method to extract diagonal asymptotics to be found in the literature 
is what we call 
the \emph{standard diagonal method}.
It consists of two steps:
(1) find a closed form or defining equation for the diagonal generating function
$G(x):= \sum_{n\in\nats} f_{a_1n,\ldots,a_dn} x^{n}$
via contour integration (in the style of \cite{Furs1967} and \cite{HaKl1971}
and summarized in \cite[Section 6.3]{Stan1999});
(2) apply univariate singularity analysis to $G(x)$ 
(in the style of \cite{FlSe} for instance) 
to compute the asymptotics.
However, as elaborated below, this method is quite limited, working well
only for main diagonals in two variables,
that is, for computing the asymptotics of $f_{nn}$.

In this article we present a new, more versatile diagonal method.
It consists of one step: apply multivariate singularity analysis
(in the style of \cite{PeWi2002}, \cite{PeWi2004}, \cite{BaPe})
directly to $F(\xvec)$ to compute the asymptotics.

\section{Limitations of the Standard Diagonal Method}\label{standard}

What is wrong with the standard diagonal method?
It certainly works well for computing main diagonal asymptotics in 
two variables, as witnessed by the following example taken from \cite{MuZa2004}. 

\begin{example}[Zigzag-free  Binary Words]
The bivariate generating function
\[
  F(x,y) = \sum_{m,n} f_{mn} x^m y^n 
     = \frac{1 + xy + x^2 y^2}{1- x- y+ xy- x^2 y^2} 
\]
counts the number of words over a binary alphabet, $\{0,1\}$ say, 
that have $m$ zeros and $n$ ones and do not contain zigzags, that is,
the subwords $010$ and $101$.
The main diagonal coefficients $f_{nn}$, then, count zigzag-free binary words 
with an equal number of zeros and ones.

To compute the asymptotics of $f_{nn}$ using 
the standard diagonal method we proceed as follows.
Since $F(x,y)$ is rational and holomorphic in a neighborhood of the origin, 
for fixed $x$ small enough $F(x/t,t)$ will be rational and holomorphic
as a function of $t$ in some annulus about $t=0$.
Thus in that annulus it can be represented by a Laurent series whose constant
term is $[t^0] F(x/t,t) = \sum_{n\ge 0} f_{nn} x^{n}$, 
the series we want.
By Cauchy's Integral and Residue Theorems we have that for some circle 
$\gamma_x$ about $t=0$
\begin{align*}
  \sum_n f_{nn} x^n 
  &= [t^0] F(x/t,t) \\
  &= \frac{1}{2\pi\mi} \int_{\gamma_x} \frac{F(x/t,t)}{t} \dif{t} \\
  &= \sum_k \Res(F(x/t,t)/t; t=s_k),
\end{align*}
where the $s_k$ are the ``small" singularities of $F(x/t,t)/t$, that is,  
the ones satisfying $\lim_{x\to 0} s_k(x) = 0$.
Since $F$ is rational, these singularities are poles and algebraic functions
of $x$, so that the residue sum, the diagonal generating function,
is also an algebraic function of $x$.

In particular,
\[
  F(x/t,t)/t =  \frac{1+ x+ x^2}{-t^2+ (1+ x- x^2)t- x},
\]
which has a single simple pole approaching zero as $x$ approaches $0$, namely, 
$s= \frac{1}{2}(1+ x- x^2- \sqrt{1- 2x- x^2- 2x^3+ x^4})$.
The residue, the diagonal generating function, is then 
\[
  G(x):= \sqrt{\frac{x^2+ x+ 1}{x^2- 3x+ 1}}.
\]

Now, the singularity of $G$ closest to the origin is $\omega:= (3-\sqrt{5})/2$,
and its reciprocal is the exponential growth order of the coefficients of $G$.
To determine the leading subexponential factor we note that 
\[
  G(x) 
  \sim  (1- \frac{x}{\omega})^{-1/2} 
     \left( 
     \sqrt{\frac{x^2+ x+ 1}{-\omega \left( x- \frac{3+ \sqrt{5}}{2} \right)}} 
     \right)_{x=\omega}
  =  (1- \frac{x}{\omega})^{-1/2}  \frac{2}{\sqrt{\sqrt{5}}}
\]
as $x\to\omega$, so that
\[
  f_{nn} = g_n 
  \sim  \omega^{-n} \frac{2}{\sqrt{\sqrt{5}}} \frac{n^{1/2-1}}{\Gamma(1/2)}
  =  \left( \frac{2}{3-\sqrt{5}} \right)^n 
     \frac{2}{\sqrt{\sqrt{5} \pi n}}
\]
by asymptotic transfer (see \cite[Chapter VI]{FlSe} for instance).
\exsym
\end{example}

However, the standard diagonal method encounters major problems off the main diagonal 
even in two variables, as illustrated in the next example adapted from \cite[Section 6.3]{Stan1999}.

\begin{example}[Lattice Paths]
The bivariate generating function
\[
  F(x,y) = \frac{1}{1- x- y- xy},
\]
whose coefficients $f_{mn}$ are called the Delannoy numbers, 
counts the number of lattice paths from $(0,0)$ to $(m,n)$ with allowable steps 
$(1,0)$, $(0,1)$, and $(1,1)$.

To compute the asymptotics of the general diagonal coefficients
$f_{an,bn}$ using the standard diagonal method,
we fix $x>0$ small enough and try to find the small poles of
\begin{align*}
  G(x)
  &:= \sum_n f_{an,bn} x^n \\
  &= [t^0] F(x^{1/a}/t^b,t^a) \\
  &= \frac{1}{2\pi\mi} \int_{\gamma_x} \frac{F(x^{1/a}/t^b,t^a)}{t} \dif{t} \\
  &= \frac{1}{2\pi\mi} \int_{\gamma_x}
     \frac{t^{b-1}}{t^{b}- x^{1/a}- t^{a+b}- x^{1/a}t^a} \dif{t}
\end{align*}
Since $t^{b}- x^{1/a}- t^{a+b}- x^{1/a}t^a$ has a zero of
multiplicity $b$ at $t=0$ when $x=0$, it follows that $F(x^{1/a}/t^b,t^a)/t$
has a single pole $s$ of order $b$ satisfying $\lim_{x\to 0} s(x) = 0$.
Thus we have that
\begin{align*}
  G(x)
  &= \Res(F(x^{1/a}/t^b,t^a)/t; t=s) \\
  &= \lim_{t\to c} \frac{1}{(b-1)!} \Dif_t^{b-1} 
     \left( (t-s)^b F(x^{1/a}/t^b,t^a)/t \right),
\end{align*}
where $\Dif_t$ is the derivative with respect to $t$.
Patiently tracing through Leibniz's rule for the iterated 
derivative of a product, we could express this limit in terms of $s$ and 
$g(t):= (t-s)^{-b}(t^{b}- x^{1/a}- t^{a+b}- x^{1/a}t^a)$,
and use it to find an algebraic equation satisfied by $G$. 
However, even with the help of Maple, this seems unlikely for general
$a$ and $b$.

Thus, already at the first step of the standard diagonal method,
it seems we are thwarted, the cause of the problem being the parameters
$a$ and $b$ occurring in the exponent: the larger they get, 
the greater the difficulty in finding a defining equation for $G$.
\exsym
\end{example}

So we should keep to main diagonal asymptotics when using the standard
diagonal method.
But even this poses problems in three or more variables,
as illustrated by the next example.

\begin{example}[Ternary Words]
The trivariate generating function
\[
  F(x,y,z) = \frac{1}{1-x-y-z}
\]
counts the the number of words $f_{lmn}$ over a ternary alphabet, 
$\{0,1,2\}$ say, that have $l$ zeros, $m$ ones, and $n$ twos
(an easy consequence of the symbolic method as described in \cite{FlSe}).
An easy combinatorial argument shows that $f_{lmn} = \binom{l+m+n}{l,m,n}$,
to which one could apply Stirling's formula and derive the asymptotics.

To compute the asymptotics of the main diagonal coefficients $f_{nnn}$ using
the standard diagonal method instead, 
we iterate the contour integration process.
First 
\begin{align*}
  \sum_n f_{m,n,n} x^m y^n 
  &= [t^0] F(x,y/t,t) \\
  &= \frac{1}{2\pi\mi} \int_{\gamma_x} \frac{F(x,y/t,t)}{t} \dif{t} \\
  &= \frac{1}{2\pi\mi} \int_{\gamma_x} 
  \frac{1}{-t^2+ (1-x)t- y} \dif{t},
\end{align*}
which has a single simple pole approaching zero as $y$ and $x$ approach $0$,
namely, $s= \frac{1}{2}(1+ x- \sqrt{1- 2x+ x^2- 4y})$.
The residue is then 
\[
  G(x,y):= \frac{1}{\sqrt{1- 2x+ x^2- 4y}}.
\]
Second, since
\[
  \lim_{x\to 0} G(x/t,t) 
   = \lim_{x\to 0} \frac{1}{\sqrt{1- 2x/t+ x^2/t^2- 4t}}
   = \frac{1}{\sqrt{1- 4t}}, 
\]
$G(x/t,t)$ is a holomorphic function of $t$ in some annulus about $t=0$. 
Thus computing 
\[
  \frac{1}{2\pi\mi} \int_{\gamma_x} \frac{G(x/t,t)}{t} \dif{t} 
  =  \frac{1}{2\pi\mi} \int_{\gamma_x} \frac{1}{t \sqrt{1- 2x/t+ x^2/t^2- 4t}}
     \dif{t}
\]
will give us the diagonal generating function.
Proceeding, we employ Maple and find that two singularities of $G(x/t,t)/t$ 
approach $0$ as $x$ approaches $0$, namely
\[
 -\frac{1}{24} p(x)+ \frac{q(x)}{p(x)}+ \frac{1}{12} 
 + \frac{1}{2}\mi\sqrt{3} 
 \left( \frac{1}{12}p(x) + \frac{q(x)}{p(x)} \right) \\
\]
and its conjugate, where 
$p(x) = (-36x+ 216x^2+ 1+ 24\sqrt{-3x^3+ 81x^4})^{1/3}$ and 
$q(x) = x-\frac{1}{24}$.
While both of these singularities are algebraic in $x$, 
Maple can not compute the sum of their residues in reasonable time.

The problem is that the residue sum, that is, 
the diagonal generating function $\sum\binom{3n}{n,n,n} x^{n}$, 
is not algebraic (see \cite[Exercise solution 6.3]{Stan1999} for instance).
While it is D-finite  ---the diagonal of any rational function 
in any finite number of variables is D-finite \cite{Lips1988}--- 
it not trivial to find a closed form or defining differential equation
for it.  
(In the present case we obtain, using the Maple package Mgfun, the defining differential equation 
$27x(x-1)f''(x)+(54x-1)f'(x)+6f(x)=0$, which has a hypergeometric solution).
Moreover, even given the defining differential equation, we still require the theory of 
univariate singularity analysis of D-finite functions to compute the asymptotics.
But this theory has not been worked out in general.
Indeed, as Philippe Flajolet has informed us, 
certain aspects of it, such as the so-called connection problem,
might not be computable!
\exsym
\end{example}

\section{A New Diagonal Method}\label{new}

To transcend the limitations of the standard diagonal method, 
to go beyond main diagonal asymptotics of bivariate generating functions,
we take a multivariate approach. When the diagonal method works, it does at least 
produce an explicit algebraic formula for the diagonal generating function, which our 
method here does not. However, it is normally much more useful to have an asymptotic 
expression for the coefficients. In order to achieve this, we can apply a multivariate 
singularity analysis directly to $F(\xvec)$.
Such an analysis has been recently developed by Baryshnikov, Pemantle, 
and Wilson in \cite{PeWi2002}, \cite{PeWi2004}, and \cite{BaPe},
the relevant parts of which we now summarize and adapt to our needs.\footnote{
As a whole, the multivariate singularity analysis developed in the
aforementioned articles  is much more general than what we present 
in this article, as it applies not only to rational functions 
but also to locally meromorphic functions. 
For a more complete summary including further examples, see the forthcoming
article \cite{PeWia}.
}

Let $\s{D}\subset\comps^d$ be the open domain of convergence of $F$ and
write $F(\xvec) = I(\xvec)/J(\xvec)$ for some $I$ and $J$ holomorphic on an
open domain $\s{D^\prime}$  containing the closure $\supsetneq \s{D}$, and relatively prime
in the ring of holomorphic functions on $\s{D^\prime}$. In all examples in this paper, the 
representation of $F$ as a quotient in fact holds on all of $\comps^d$, but the extra generality is 
sometimes useful in applications.
Let \s{V}\ be the complex variety $\set{\xvec\in\comps^d}{J(\xvec)=0}$.

A \hl{critical point} of $F$ for $\nvec\in(\pnats)^d$ is a solution of   
\begin{align*}
  J(\xvec) 
     &= 0 \\
  n_d x_i J_i(\xvec) 
     &= n_i x_d J_d(\xvec) \quad(i<d).
\end{align*}
Here $n_i$ denotes component $i$ of \nvec, and $J_i$ denotes the 
partial derivative of $J$ with respect to component $i$ of the domain of $J$,
conventions we adopt for all vectors and functions throughout 
and which should cause no confusion in context.
Let \crit{\nvec}\ denote the set of all critical points of $F$ for $\nvec$. For generic directions 
$\nvec$, this set is finite, being a zero-dimensional complex variety. The main situation in 
which \crit{\nvec}\  is infinite  occurs when $J$ defines a binomial variety
$\{ \xvec \mid \xvec^\mathbf{a} - \xvec^\mathbf{b} \}$, in which case \crit{\nvec}\  is empty for 
all but one direction and uncountable otherwise. Such examples can be analysed by a variant of the 
methods shown here.

A \hl{contributing point} of $F$ for $\nvec$ is a critical point that 
influences the asymptotics of the coefficients of $F$ in the direction of \nvec.
Let \contrib{\nvec}\ denote the set of all such points.
While \contrib{\nvec}\ is ill-defined here, 
its functional role will become clear from the next two theorems.

\begin{theorem}\label{contrib}
If \crit{\nvec}\ is finite, then
\begin{itemize}
\item \crit{\nvec}\ contains exactly one point, call it \cvec, that lies in 
  the positive orthant of $\reals^d$, and $\cvec\in \contrib{\nvec}$;
\item all other members of \contrib{\nvec} must lie on the same torus as \cvec;
\item in the case where $J = 1 - P$ for some aperiodic power series 
  with nonnegative coefficients $P$, $\contrib{\nvec} = \{\cvec\}$.
\end{itemize}
\end{theorem}

Here, the \hl{torus} of a point $\cvec\in\comps^d$ is the set 
$\set{\xvec\in\comps^d}{\forall i\le d \; |x_i| = |c_i|}$, 
and a power series is \hl{aperiodic} if the $\ints$-span of its
monomial vectors is all of $\ints^d$.

A point $\cvec\in\s{V}$ is a \hl{smooth point} if \s{V}\ is a smooth 
complex manifold in a neighborhood of \cvec, or equivalently, 
if $J_i(\cvec) \not= 0$ for some $i$
(see \cite[page 363]{BrKn1986} for instance).
For simplicity of presentation, we deal only with smooth points in this article. 
This is the generic case, 
although interesting examples are not always generic. For more on the case of non-smooth 
points, see \cite{PeWi2004}.

\begin{theorem}\label{asymptotics}
Let $\nvec = (a_1 n,\ldots,a_d n)$ for some $a_1,\ldots,a_d\in\pnats$.
If \contrib{\nvec}\ consists of a single smooth point $\cvec$ such that 
$c_d J_d(\cvec) \not= 0$ and
$c_d$ is a simple zero of $x_d \mapsto J(c_1,\ldots,c_{d-1},x_d)$, then
\[
  f_\nvec 
  =  \cvec^{-\nvec} \left( 
     \sum_{j = 0}^{j^\prime} b_j n_d^{(1-d-j)/2} 
     + O\left( n_d^{-d-j^\prime} \right)
     \right)
\]
for some $b_j\in\comps$ as $n\to\infty$.
If \contrib{\nvec}\ consists of a finite set of smooth points
$\cvec_1,\ldots,\cvec_r$ each satisfying the hypotheses above,
then
\[
  f_\nvec 
  =  \sum_i \cvec_i^{-\nvec} \left(
     \sum_{j = 0}^{j'} b_{ij} n_d^{(1-d-j)/2} 
     + O\left( n_d^{-d-j'} \right)
     \right)
\]
as $n\to\infty$.
\end{theorem}

\begin{theorem}\label{hessian}
In Theorem~\ref{asymptotics}, 
\[
  b_0 = \frac{I(\cvec)}{-c_d J_d(\cvec) \sqrt{(2\pi)^{d-1} h(J,\cvec)}},
\]
assuming $h(J,\cvec) \not= 0$.
Here $h(J,\cvec)$ is the determinant of the matrix \hmatrix, whose entries are
\begin{align*}
  \hmatrix_{lm}
  &= \frac{c_l c_m}{c_d^2 J_d^2}
     \left( J_m J_l+ c_d 
     (J_d J_{lm}- J_m J_{ld}- J_l J_{dm}+ \frac{J_l J_m}{J_d} J_{dd})
     \right) \Big|_{\xvec = \cvec}; \\
  \hmatrix_{ll}
  &= \frac{c_l J_l}{c_d J_d} + \frac{c_l^2}{c_d^2 J_d^2}
     \left( J_l^2+ c_d 
     (J_d J_{ll}- 2 J_l J_{ld}+ \frac{J_l^2}{J_d} J_{dd})
     \right) \Big|_{\xvec = \cvec}, 
\end{align*}
where $l,m < d$ and $l\not= m$.
\end{theorem}

\begin{proof}
The formula for $b_0$ comes from \cite[Theorem 3.5]{PeWi2002}.
(Please note the typo therein: the $z_d H_d$ should be a $-z_d H_d$.)
We prove the formula for \hmatrix\ only.

For easy reading, let $\xhat = (x_1,\ldots,x_{d-1})$.
Since $J$ is holomorphic and $J_d(\cvec) \not= 0$, 
there exists a holomorphic function $g$ guaranteed by the
Implicit Function Theorem such that in some open ball around \cvec,
$J(\xhat,g(\xhat)) = 0$.
As shown in \cite{PeWi2002}, \hmatrix\ is the Hessian
(matrix of second partial derivatives) of
\[
\ftilde(\thetahat):= 
     +\log g(c_1\me^{\mi\theta_1},\ldots,c_{d-1}\me^{\mi\theta_{d-1}})
     + \mi \sum_{j=1}^{d-1} \frac{n_j}{n_d}\theta_j - \log g(\cvec)
\]
evaluated at $\theta = \zerovec$.

Now, from the Implicit Function Theorem we also get
\begin{align*}
  \Dif_l g(\xhat) 
     &= -\frac{J_l(\xhat,g(\xhat))}{J_d(\xhat,g(\xhat))} 
        \quad (l<d);\\
  \Dif_m J_l (\xhat,g(\xhat)) 
     &= J_{lm} + J_{ld} \Dif_m g 
      = J_{lm} - J_{ld} \frac{J_m}{J_d}
        \quad (l\le d, m<d).
\end{align*}
Thus, up to an additive constant,
$\Dif_l \ftilde = -\mi c_l \me^{i\theta_l} g^{-1} J_l/J_d$,
and for $l,m < d$ with $l\not= m$ we have 
\begin{align*}
  \Dif_{lm} \ftilde \Big|_{\thetahat = \zerovec}
     &= c_l c_m 
        \left(-g^{-2} \Dif_m g \frac{J_l}{J_d}+ 
        g^{-1} \frac{J_d \Dif_m J_l- J_l \Dif_m J_d}{J_d^2} 
        \right) \Big|_{\xvec = \cvec}; \\
  \Dif_{ll} \ftilde \Big|_{\thetahat = \zerovec}
     &= c_l g^{-1} \frac{J_l}{J_d}+  c_l^2 
        \left(-g^{-2} \Dif_l g \frac{J_l}{J_d}+ 
        g^{-1} \frac{J_d \Dif_l J_l- J_l \Dif_l J_d}{J_d^2} 
        \right) \Big|_{\xvec = \cvec}. \\
\end{align*}
The result then follows by plugging in the Implicit Function Theorem equations
above and simplifying.
\end{proof}

Several of our examples involve main diagonal asymptotics for symmetric $J$.
In this case $h(J,\cvec)$ simplifies greatly.

\begin{proposition}\label{symmetric}
If \crit{\nvec}\ is finite, $\nvec = (n,\ldots,n)$,
$J(\xvec)$ is symmetric in \xvec, and $J_d(\cvec) \not= 0$, 
where $\cvec$ is the contributing point 
that lies in the positive orthant of $\reals^d$, 
then $\cvec = (c,\ldots,c)$ for some positive real $c$, and
\[
  h(J,\cvec) = d \left( 1 + \frac{c}{J_d}(J_{dd} - J_{d1})\Big|_{\xvec = \cvec}
     \right)^{d-1}.
\]
\end{proposition}

\begin{proof}
By Theorem~\ref{contrib}, \crit{\nvec}\ 
contains exactly one point \cvec\ that lies in the positive orthant of 
$\reals^d$.
By the symmetry of the critical equations $J=0$ and $x_i J_i = x_d J_d$
$(i<d)$ induced by the symmetry of \nvec\ and $J$,
any permutation of \cvec's coordinates is also a critical point
lying in the positive orthant of $\reals^d$.
Since there can be only one such point, $\cvec = (c,\ldots,c)$ for some
positive real $c$.

Since $J$ is symmetric and $\cvec = (c,\ldots,c)$, 
all first partial derivatives of $J$ are equal at \cvec,
all mixed second partial derivatives of $J$ are equal at \cvec, and 
all diagonal second partial derivatives of $J$ are equal at \cvec.
Thus the entries of the matrix of Theorem~\ref{hessian} simplify to
\begin{align*}
  \Dif_{lm} \ftilde \Big|_{\thetahat = \zerovec}
     &= 1 + \frac{c}{J_d}(J_{dd} - J_{d1}) \Big|_{\xvec = \cvec} \\
  \Dif_{ll} \ftilde \Big|_{\thetahat = \zerovec}
     &= 2 \Dif_{lm} \ftilde \Big|_{\thetahat = \zerovec},
\end{align*}
where $l,m < d$ and $l\not= m$, and the determinant simplifies to
$h(J,\cvec) = a^{d-1}d$,
where $a = 1 + \frac{c}{J_d}(J_{dd} - J_{d1}) \Big|_{\xvec = \cvec}$.
\end{proof}

Our new diagonal method is simply the application of 
Theorems~\ref{contrib}, \ref{asymptotics}, and \ref{hessian} directly to $F$.
We illustrate this now with several examples.

\begin{example}[Zigzag-free Binary Words]
Consider again the bivariate generating function 
\[
  F(x,y) = \sum_{m,n} f_{mn} x^m y^n 
  = \frac{1 + xy + x^2 y^2}{1- x - y + xy - x^2 y^2}
\]
of Section~\ref{standard}.
To compute the asymptotics of $f_{nn}$ using the new diagonal method,
let $\nvec = (n,n)$ and $J = 1- x- y+ xy- x^2 y^2$.
Then \crit{\nvec}, the solution set of $J = 0$ and $nx J_x = ny J_y$, comprises
$(1/\phi, 1/\phi)$, $(-\phi,-\phi)$, $(1\pm \mi\sqrt{3})/2$, and 
$(1\pm \mi\sqrt{3})/2$, where $\phi = (1+\sqrt{5})/2$.
By Theorem~\ref{contrib}, $\cvec:= (1/\phi, 1/\phi)$ is a contributing point
and the only such point, since none of the other critical points lie on the
same torus as \cvec.
It is also a smooth point. 
By Proposition~\ref{symmetric}, $h(J,\cvec) = 4/(-5+ 3\sqrt{5})$, and
by Theorems~\ref{asymptotics} and \ref{hessian},
\[
  f_{nn} \sim \phi^{2n} \frac{2}{\sqrt{\sqrt{5} \pi n}}.
\]
This agrees with our answer in Section~\ref{standard} 
since $\phi^2 = 2/(3- \sqrt{5})$.
\exsym
\end{example}

\begin{example}[Lattice Paths]
Consider again the bivariate generating function
\[
  F(x,y) = \frac{1}{1- x- y- xy},
\]
of Section~\ref{standard}.
To compute the asymptotics of $f_{an,bn}$ using the new diagonal method,
let $\nvec = (an,bn)$ and $J = 1- x- y- xy$.
Then \crit{\nvec}, the solution set of $J = 0$ and $bnx J_x = any J_y$,
comprises $(\frac{-a\pm L}{b},\frac{-b\pm L}{a})$, where $L=\sqrt{a^2+b^2}$.
By Theorem~\ref{contrib}, $\cvec:= (\frac{-a+ L}{b},\frac{-b+ L}{a})$
is the only contributing point.
It is also a smooth point.
By Theorem~\ref{hessian}, $h(J,\cvec) = \frac{-2(b- L)a(a^2+b^2- aL)}
{(a- L)^2 (a- b+ L)^2)}$, and by Theorem~\ref{asymptotics},
\[
  f_{an,bn} 
  \sim \left( \frac{-a+ L}{b} \right)^{-an} 
        \left( \frac{-b+ L}{a} \right)^{-bn} 
        \sqrt{\frac{ab}{L (a+ b- L)^2 2\pi n}}.  
\]
Here, where the standard diagonal method failed, 
the new diagonal method allowed us to compute the asymptotics quite easily
(with the help of Maple).
This success comes from a general attribute of the new diagonal method: 
the parameters $a$ and $b$ appear as factors instead of as powers
in the equations we need to solve. 
\exsym
\end{example}

\begin{example}[Ternary Words]
Consider again the trivariate generating function 
\[
  F(x,y,z) = \frac{1}{1- x- y- z}
\]
of Section~\ref{standard}.
To compute the asymptotics of $f_{an,bn,cn}$,
let $\nvec = (an,bn,cn)$ and $J = 1- x- y- z$.
Then \crit{\nvec}, the solution set of $J = 0$, $cnx J_x = anz J_z$, and
$cny J_y = bnz J_z$ comprises $\cvec:= (a+b+c)^{-1} (a,b,c)$.
By Theorem~\ref{contrib}, $\cvec$ is the only contributing point. 
It is also a smooth point.
By Proposition~\ref{symmetric}, $h(J,\cvec) = ab(a+b+c)/ c^3$, 
and by Theorems~\ref{asymptotics} and \ref{hessian},
\[ 
  f_{an,bn,cn} \sim 
     \frac{(a+b+c)^{(a+b+c)n}}{a^{an}b^{bn}c^{cn}}
     \sqrt\frac{a+b+c}{abc} \frac{1}{2\pi n},
\]
in agreement with Stirling's formula.

Here, where the standard diagonal method failed because of the nonalgebraic
character of the diagonal generating function, 
the new diagonal method gave the answer easily, even off the main diagonal.
\exsym
\end{example}

Rather than multiply examples unnecessarily, we finish with an example involving an 
arbitrary number of variables. Many more examples of multivariate asymptotics using 
similar methods can be found in \cite{PeWia}.

\begin{example}[Alignments]
In a problem from computational biology, one is
given a finite set of strings of various lengths over some
finite alphabet and needs to count the number of ways of inserting blanks
into these strings to make them all of equal length.
The strings represent genetic sequences, and inserting blanks represents
aligning genetically similar segments among the sequences.

More specifically and abstractly, a \emph{$(d,r_1,\ldots,r_d)$-alignment}
is a table of $d$ $\{0,1\}$-strings, one in each row, such that 
the sum of the entries of row $j$ equals $r_j$ $(1\le j\le d)$ and 
no column contains all zeros.  
(The ones correspond to letters of the genetic sequence 
and the zeros to blanks.)
As shown in \cite{GHOW1990}, 
the generating function for the number of $(d,\cdot)$-alignments is
\[
  F(x_1,\ldots,x_d) = \frac{1}{2- \prod_{i=1}^d (1+x_i)}.
\]

To compute the asymptotics of $f_{n,\ldots,n}$   
(the case where all rows have the same sum;
the genetic sequences all have the same length) using the new diagonal method,
let $\nvec = (n,\ldots,n)$, $I = 1/2$ and $J = 1- (1/2)\prod_{i=1}^d (1+x_i)$.
Then \crit{\nvec}, the solution set of $J = 0$ and $nx_i J_i = nx_d J_d$ 
($i<d$) comprises $(2^{1/d}\omega_i-1,\ldots,2^{1/d}\omega_i-1)$,
where the $\omega_i$ ($i<d$) are the $d$th roots of unity.
By Theorem~\ref{contrib}, $\cvec:= (2^{1/d}-1,\ldots,2^{1/d}-1)$
is the only contributing point. 
It is also a smooth point.
By Proposition~\ref{symmetric}, $h(J,\cvec) = d 2^{(1-d)/d}$, and
by Theorems~\ref{asymptotics} and \ref{hessian},
\[
  f_{n\ldots n} \sim (2^{1/d}-1)^{-dn} 
  \frac{1}{(2^{1/d}-1) 2^{(d^2-1)/2d} \sqrt{d (\pi n)^{d-1}}}, 
\]
in agreement with \cite{GHOW1990}. 
(Please note the typo therein: the $2^{(k^2-1)/2k}$ in the formula on page 139
should be a $2^{(1-k^2)/2k}$).

A more important case for biological applications is to count alignments
whose minimum block size is bounded below by a constant.
A \emph{block} of size $b$ in a $(d,\cdot)$-alignment is a $d\times b$ 
submatrix of the alignment with contiguous columns,
none of which contain a zero.
As shown in \cite{ReTa}, the generating function for $(d,\cdot)$-alignments with
blocks of size at least $b$ is
\begin{align*}
  G(x_1,\ldots,x_d) 
  &= \frac{1- t+ t^b}{1+ (1-p)(1-t+ t^b)- t^2+ t^{b+1}} \\
  &= \frac{1- t+ t^b}{1- t- (p- 1- t)(1- t+ t^b)}
\end{align*}
where $p= \prod_{i=1}^d (1+x_i)$ and $t= \prod_{i=1}^d x_i$.
Note that when $b=1$, $G = F$, as expected.

To compute the asymptotics of $g_{n,\ldots,n}$
(which to our knowledge has only been done up till now for $d=2$ 
via the standard diagonal method \cite{ReTa}) 
let $\nvec= (n,\ldots,n)$ and $J = 1- t- (p- 1- t)(1- t+ t^b)$.
Then \crit{\nvec}\ comprises all vectors $x\onevec$ such that
$j(x):= J(x\onevec) = 0$ 
(as seen from the equations $nx_i J_i = nx_d J_d$).
Let $\cvec:= c\onevec$ be the critical point in the positive orthant 
of $\reals^d$ (so $c$ is the unique positive solution of $j(x)=0$) 
and note that  $c < 1$, because $j(x) < 0$ for $x\ge 1$.
By Theorem~\ref{contrib}, $\cvec\in\contrib{\nvec}$ and 
all other points of \contrib{\nvec}\ lie on the same torus as \cvec.
Since  
\[
  G(x_1,\ldots,x_d) 
  = \frac{1+ \frac{t^b}{1-t}}{1-(p-1-t)(1+ \frac{t^b}{1-t})},
\] 
on the open polydisk $\set{\xvec\in\comps^d}{\forall i \; |x_i|< 1}$,
since this polydisk contains the torus of \cvec, and since 
$(p-1-t)(1+ \frac{t^b}{1-t})$ is an aperiodic power series with nonnegative 
coefficients, we can apply Theorem~\ref{contrib} to conclude that \cvec\ 
is the only contributing point.
It is also a smooth point, as is seen most easily from the form $J = 1-(p-1-t)(1+ \frac{t^b}{1-t})$.
Thus by Theorem~\ref{asymptotics},
\[
  g_{n\ldots n} \sim c^{-dn} b_0 n^{(1-d)/2},
\]
where $b_0$ is the constant given in Theorem~\ref{hessian}.
Note that both $c$ and $b_0$ depend on $b$.
\exsym
\end{example}

\bibliographystyle{amsalpha}    
\bibliography{mgf}
\end{document}